# Optimizing Return and Secure Disposal of Prescription Opioids to Reduce the Diversion to Secondary Users and Black Market

**Md Mahmudul Hasan**
Ph.D.
Dept. of Pharmaceutical Outcomes and Policy
College of Pharmacy
Dept. of Information Systems and Operations Management
College of Business
University of Florida
1225 Center Drive, Gainesville, FL 32606, USA
Email: hasan.mdmahmudul@ufl.edu

**Tasnim Ibn Faiz**
Ph.D.
Dept. of Mechanical and Industrial Engineering
College of Engineering
Northeastern University
360 Huntington Avenue, Boston, MA 02135, USA
Email: faiz.t@northeastern.edu

**Alicia Sasser Modestino**
Ph.D.
College of Social Science and Humanities
Center for Health Policy and Healthcare Research
Northeastern University
360 Huntington Avenue, Boston, MA 02135, USA
Email: a.modestino@northeastern.edu

**Gary J. Young**
J.D., Ph.D.
D'Amore-McKim School of Business
Bouve College of Health Sciences
Center for Health Policy and Healthcare Research
Northeastern University
360 Huntington Avenue, Boston, MA 02135, USA
Email: ga.young@northeastern.edu

**Md. Noor-E-Alam**
Ph.D.
Dept. of Mechanical and Industrial Engineering
College of Engineering
Center for Health Policy and Healthcare Research
Northeastern University
360 Huntington Avenue, Boston, MA 02135, USA
Email: mnalam@northeastern.edu

**Corresponding author**
**Md. Noor-E-Alam**
Ph.D.
Dept. of Mechanical and Industrial Engineering
College of Engineering
Center for Health Policy and Healthcare Research
Northeastern University
360 Huntington Avenue
Boston, MA 02135, USA
(617) 373-2275
Email: mnalam@northeastern.edu

**Financial disclosure and funding source:** The research was supported by a grant from Northeastern University's Global Resilience Institute (GRI).

# Abstract

Opioid Use Disorder (OUD) has reached an epidemic level in the US. Diversion of unused prescription opioids to secondary users and black market significantly contributes to the abuse and misuse of these highly addictive drugs, leading to the increased risk of OUD and accidental opioid overdose within communities. Hence, it is critical to design effective strategies to reduce the non-medical use of opioids that can occur via diversion at the patient level. In this paper, we aim to address this critical public health problem by designing strategies for the return and safe disposal of unused prescription opioids. We propose a data-driven optimization framework to determine the optimal incentive disbursement plans and locations of easily accessible opioid disposal kiosks to motivate prescription opioid users of diverse profiles in returning their unused opioids. We develop a Mixed-Integer Non-Linear Programming (MINLP) model to solve the decision problem, followed by a reformulation scheme using Benders Decomposition that results in a computationally efficient solution. We present a case study to show the benefits and usability of the model using a dataset created from Massachusetts All Payer Claims Data (MA APCD). Our proposed model allows the policymakers to estimate and include a penalty cost considering the economic and healthcare burden associated with prescription opioid diversion. Our numerical experiments demonstrate the model's ability and usefulness in determining optimal locations of opioid disposal kiosks and incentive disbursement plans for maximizing the disposal of unused opioids. The proposed optimization framework offers various trade-off strategies that can help government agencies design pragmatic policies for reducing the diversion of unused prescription opioids.

## 1. Introduction

In this manuscript, we report results from a study we conducted for reducing the non-medical use of prescription opioids that can occur via diversion of these drugs to secondary users and the black market. Patients frequently are prescribed too much opioid medication and do not consume all the pills they are prescribed [1-4]. In addition, they do not dispose of the unused drugs properly and facilitate the deliberate or inadvertent diversion of these drugs [5, 6]. Although "take-back" days have helped to collect excess opioid medications, these initiatives are not always well publicized, are limited in duration (e.g., two days per year) or are located at inconvenient sites (e.g., police stations). We sought to address this issue by proposing a unique scheme that includes a two-fold objective: (i) to design incentive allocation policies for opioid users (opioid users and patients are used interchangeably throughout the manuscript) with different characteristics, and (ii) to locate opioid disposal kiosks at a site convenient to patients to ensure safe and secure disposal of unused opioids.

To assess the feasibility of this proposed scheme and analyze pertinent policies, we developed a decision-making framework using optimization methods. In particular, we formulated a *Mixed-Integer Non-Linear Programming* (MINLP) model that can generate optimal decisions pertaining to the location of opioid disposal kiosks and the amount of a financial incentive to the opioid users for returning the unused drugs. The proposed model demonstrated ability for making various tradeoffs among these two decisions and also for providing an estimate of the overall cost needed to run the return campaign. To our knowledge, this is the first study that develops a data-driven optimization framework to address the crisis of prescription opioid diversion that generally occurs at the patient-level. In our study, we constructed patients' profiles based on the risk factors associated with OUD as well as based on the factors that influence the patients' opioid return tendency. We also included patients' specific minimum reservation incentives in our model to assess the influence of patients' willingness to return their unused opioids on various trade-off strategies. In addition, our modeling framework provides a unique opportunity for public health policy makers to consider the penalty cost that we specifically included considering the societal, environmental, health, and economic burden associated with the misuse and abuse of diverted prescription opioids to the community. Our proposed decision-making framework, following further validation, can help government agencies design implementable policies to reduce the diversion of unused prescription opioids by coordinating efforts among different stakeholders (e.g., pharmaceutical manufacturing companies, local pharmacies, community clinics, and other third-party organizations).

## 2. Background and Literature Review

The misuse of prescription opioids and opioid use disorder (OUD) have created a global public health crisis, affecting the lives of individuals in many countries including North America, West and North Africa, Near and Middle East, Asia, and Western and Central Europe [7]. Across the globe, an estimated 62 million people aged 15 to 64 reported non-medical use of opioids (i.e., prescription opioids, heroin, and synthetic opioids) in 2019 [7], and over 115,000 died of an opioid overdose in 2017 [8]. This public health problem has become even more severe in the US in recent years. In 2020, an estimated 2.7 million Americans aged 12 and older suffered from OUD [9]. OUD significantly contributes to overdoses [10], and nearly 70,000 individuals have lost their lives due to opioid overdose in 2020 [11]. Alongside this substantial death toll, the economic burden of OUD in the US is also overwhelming with an estimated annual cost exceeding $786 billion to the society in 2018 [12]. In response, many regulatory actions and evidence-based strategies have been implemented both at federal and state levels to reduce drug trafficking, inappropriate prescribing, and illegal dispensing of opioids [13-15]. Yet, the prevalence of opioid misuse and related overdose deaths have increased substantially since 1999 [1]. ED visits due to accidental ingestion and prescription-opioid-abuse related treatment admissions have also reportedly increased over the last decades [16]. Some studies have reported that roughly 80% of injection drug users have abused prescription medications, suggesting that prescriptions may sometimes serve as a gateway for injecting heroin and other drugs [16-19].

Prior studies reported that opioids may be prescribed in an inappropriately higher number and outside of medical need [16]. In 2012, 259 million opioid prescriptions were dispensed—effectively enough for every American to have a bottle of pills in their medicine cabinet [20]. Although a legitimate need exists for opioids in certain clinical conditions (e.g., chronic pain and surgical procedures), prior studies have indicated that a significant number of patients reported underutilization of these prescribed opioids with an intent to keep the unused medications at home [2-4]. The National Community Pharmacist Association (NCPA) also reported that up to 40% of prescription opioids are not completely used [1, 21]. These medications are likely to be kept inside the home with insecure storage, potentially posing a significant risk of diversion and accidental poisoning [1].

Such diversion of prescription opioids that occur at the patient-level is a pressing concern [1]. The results from the 2020 National Survey on Drug Use and Health (NSDUH) indicated that an estimated 9.3 million individuals aged 12 or older misused prescription pain medications in the past year [9]. Approximately 53% of these individuals obtained the medications from friends or families either for free, or for money, or took them without permission [22, 23]. Thus, ensuring proper disposal of unused prescription opioids to reduce the prevalence of misuse and minimize the risk of deliberate or inadvertent diversion is an important public health issue. As such, it is critical that steps be taken to make communities resilient against this epidemic of opioid addiction.

In this vein, a potential solution could be to provide patients with information on proper disposal at the point of dispensing the medication, to encourage disposal when the prescription is no longer being used, and to provide a financial incentive to return unused opioid medications at a site convenient to the patients. There is evidence suggesting that providing patients with some form of financial incentive and placing a disposal kiosk in a location where they frequently visit would encourage them to return their unused opioids [1, 24]. However, little is known about the factors affecting patient motivation and compliance as well as the setting in which the intervention should be delivered. To date, there has been little in-depth research investigating the impact of providing optimal incentives to legitimate opioid users in disposing and/or returning the unwanted drugs. Little is known about the complexity originated from several factors such as patients' awareness to ongoing opioid addiction epidemic, effect of financial incentives, and attitude towards returning the unused drugs that directly affect their willingness to participate in the opioid disposal campaign.

Prior studies on the design and effectiveness of Reverse Supply Chain (RSC) for pharmaceutical products indicate that the recovery of unwanted medications is complicated [25, 26]. This is mainly because manufacturers often lack information about the required budget for the collection and disposal of drugs, available quantity of unwanted medications, and consumers' willingness to return the drugs. Despite the fact that consumers play an important role for ensuring successful return of these medications, consumer behavior has not been considered in most of the prior studies [25]. In addition, none of the state regulated drug stewardship programs in the US have any actionable plans for offering some form of incentives to encourage and motivate consumers for returning unused drugs.

Several empirical studies on reverse logistics reported a *positive impact of providing incentives* to users for product recovery in re-manufacturing [27-32], and waste recycling settings [31, 33-37]. Some studies on RSC included hospitals and pharmacies as consumers and offered an incentive to them in return of the unused medications [25, 38, 39]. These studies also considered resale value of returned medications that are unexpired, resulting in revenue in the RSC for the manufacturers. However, this type of structure is not suitable for a program that focuses on the return of unused opioids given the fact that storage of unused opioids as well as their diversion and misuse occur at the patient level. In addition, considering the risk of misuse and overdose associated with the unused opioids, the potential societal, environmental, and health benefits from ensuring safe and secure return and disposal of these drugs in many ways outweigh the revenue from reselling them in the secondary market.

Some practical evidence on the collection of unused medications exists in the US. Following the adoption of the Secure and Responsible Drug Disposal Act by Congress in 2010, the Drug Enforcement Administration (DEA) in collaboration with local law enforcement and community groups launched nationwide disposal efforts at biannual Drug Take-Back Days for collecting potentially dangerous unused, unwanted and expired medications [16]. Several states including New York, California, Washington, Illinois, Vermont, and Massachusetts have passed pharmaceutical stewardship laws that require pharmaceutical manufacturers to facilitate and fund drug stewardship programs for collecting unused prescription and non-prescription drugs [40]. In contrast to the laws passed in other states, the related laws in Massachusetts and Vermont were included in broader legislation--to address state's substance use and abuse crisis. The Acts passed in these two states established drug stewardship programs that particularly apply to brand and generic opioids including benzodiazepines [40].

For example, in the state of Massachusetts, Inmar Intelligence, a third-party organization submitted a proposal to the Massachusetts Department of Public Health (MDPH) in 2020 for implementing a drug stewardship program that will be fully financed by the drug manufacturers [41]. Such a stewardship program will offer prescription medications return services via prepaid mail-back envelops and establish drug disposal kiosks in authorized collection facilities such as retail pharmacies and law enforcement agencies to facilitate the return of unused, and expired medications. However, to the best of our knowledge, the effectiveness of this program has not yet been evaluated.

In this study, we focused on the provision of incentives to multiple categories of prescription opioid users for the proper disposal of their excess drugs. We also sought to determine the optimal locations of disposal kiosks to ensure safe and secure return of unused opioids. The previously reported empirical findings [1, 24] regarding the positive impact of incentives and disposal kiosks-placement on patients' willingness to return unused opioids served as a motivation for our study. In particular, we formulated an optimization-based joint incentive and disposal-kiosk-location planning framework for determining a sufficiently motivating monetary incentive level and optimal locations of opioid disposal kiosks. Finally, we conducted a sensitivity analysis to systematically simulate, quantify and understand the implications of the proposed decision-making framework.

## 3. Problem statement

For purposes of our study, we assume a setting whereby a pharmaceutical company is legally obligated under a pharmaceutical stewardship law [40, 42] to finance drug stewardship program for the collection of unused opioids. In the event of failing to do so, the company will face some penalty imposed by the government. The pharmaceutical company contracts with a third-party stewardship organization (e.g., similar to Inmar Intelligence) to administer the campaign for collecting unused opioids at different zones. The zones could be defined as major cities where patients will return their unused opioids. As per the pharmaceutical stewardship law, the pharmaceutical companies or the stewardship organization cannot charge any point-of-sale, point-of-collection, processing fees or other drug cost to individual consumers/patients. As such, the initiative is not expected to increase the price of opioids.

We formulated the MILP model in consideration of the third-party stewardship organization's problem. Contracted by the pharmaceutical company, the third-party organization is required to identify the optimal number of disposal kiosks and amount of incentive that will be needed to ensure the return of a certain percentage of total unused opioids. In addition, the third-party organization will also need to provide an estimate of the total cost to administer this campaign. To that end, the third-party organization will install opioid disposal kiosks at a set of potential locations and offer an incentive to those using opioids who reside within certain zones to motivate them for returning the unused opioids to a designated disposal kiosk. The pharmacies that surround the patients' zone serve as the potential sites for the third-party organization to install disposal kiosks. The organization will also be aware of certain factors that potentially influence the patients' willingness to return unused opioids. Such factors include patients' age, gender, and also the type of unused opioids that patients' will potentially return. The third-party organization will categorize the patients' based on the profiles that will consist of the above-mentioned patient specific factors. Patients with different profiles will differ as to the minimum reservation incentive that will define their willingness/motivation to return unused opioids.

We begin with the third-party organization's problem of identifying the optimal number of disposal kiosks and optimal amount of incentive to offer in order to minimize the total cost. This total cost consists of the fixed cost for opening disposal kiosks, the cost of providing incentives to opioid users, and the penalty cost. While the first two cost components are self-interpretable, the penalty cost requires some explanations. We introduced the penalty cost to constrain the model to ensure sufficient return of the unused opioids commensurate with the quantity of opioids that are set as a target for the organization at the beginning of the campaign. This cost is in the form of a penalty that can be imposed by the state government on the pharmaceutical company due to the inability of returning the unused opioids potentially available in a certain zone. As noted, the opioids that are not consumed and not disposed of properly can be diverted to secondary users, potentially increasing the economic and healthcare burden in terms of accidental overdose, increased hospitalization and ER visits, and poor functional status at an individual level. Thus, the penalty cost resulting from the failure of ensuring the return of all the available unused opioids is an important consideration in the proposed framework. The basic assumptions of this modeling framework are given below:

1. Patients in accordance with their profiles have minimum reservation incentive levels, and they will return their unused opioids if the incentive is equal to or greater than the sum of that minimum reservation incentive and the transportation cost that they will incur when they dispose of the unused opioids at a disposal kiosk. This assumption implies that the opioid users will bear the travel cost to disposal kiosks. This is reasonable when the disposal kiosks are placed within the neighborhood zones of the resident.

2. The profiles of users and related minimum reservation incentives are known a priori [27, 28, 33, 43]. This assumption is reasonable given empirical findings from previous studies [44, 45] that have reported an increased risk of abusing prescription opioids in relation to patients' certain demographic characteristics such as gender (i.e., male) and age (i.e., younger individuals). Previous studies have also reported that individuals who were prescribed short-acting opioids are less likely to return unused opioids [1] indicating that patients with different characteristics will have varying response/motivation to return unused opioids.

3. An increase in the minimum reservation incentive level will motivate opioid users to travel further for returning the opioids. This is reasonable from the behavioral and economic standpoint as a larger monetary incentive serves as a driving force for motivating individuals to accomplish a task.

4. The minimum reservation incentive for different patients and the maximum distance that patients will be willing to travel is irrelative to the different patients' zones. This assumption is included to ensure that the proposed incentivization scheme can be uniformly implemented and thus is both practical and socially acceptable.

5. Given a set of available disposal kiosks, users of opioids will only be willing to return unused opioids to the nearest disposal kiosks. This assumption is reasonable in practice given that convenience in returning the opioids is usually considered as an important determinant of an individual's motivation to return unused opioids.

## 4. Mathematical model formulation

We presented the mixed-integer mathematical programming model proposed in this study in order to formulate the third-party organization's problem pertinent to the campaign of returning unused opioids under discussion. We listed the basic model parameters, and decision variables in Table 1. The notations for other variables and parameters will be introduced and defined appropriately in the following subsections.

**Table 1. Basic model parameters and decision variables**

| | |
|---|---|
| **Index sets** | |
| $I$ | Set of candidate pharmacy locations for placing an opioid disposal kiosk; $I = \{1, ..., m\}$ |
| $J$ | Set of opioid users' zones; $J = \{1, ..., n\}$ |
| $P$ | Set of profiles of opioid users; $P = \{1, ..., p\}$ |
| **Parameters** | |
| $Q_{jp}$ | Number of unused opioids available to opioid users with profile $p$ at zone $j$. |
| $\theta$ | A certain percentage of available unused opioids that the third-party organization intend to return through the campaign |
| $f_i$ | Fixed cost for opening and operating a disposal kiosk at a pharmacy located at $i$ |
| $k_i$ | Fixed capacity of the disposal kiosks |
| $c_{ij}$ | Cost of traveling from zone $j$ to a disposal kiosk placed in a pharmacy located at $i$ |
| $A_{ijp}$ | 1 if travelling distance from a zone $j$ to disposal kiosk placed at $i$ is less than the maximum distance that a user with profile $p$ will be willing to travel at given minimum reservation incentive level, or 0 otherwise |
| $m$ | Penalty cost associated with unreturned opioids |
| **Decision variables** | |
| $Y_i$ | Binary decision, 1 if a pharmacy in location $i$ is selected to open a disposal kiosk, or 0 otherwise |
| $R_{ijp}$ | Incentive paid to an opioid user with profile $p$ when he/she returns unused opioids per prescription from zone $j$ to a disposal kiosk placed at location $i$ |
| $w_{ijp}$ | Binary decision, 1 if an opioid user with profile $p$ in zone $j$ is assigned to return unused opioids to a disposal kiosk located at $i$, or 0 otherwise |
| $x_{ijp}$ | Quantity of unused opioids that are returned by users with profile $p$ from zone $j$ to a disposal kiosk located at $i$ |
| $b_{jp}$ | Unused opioids that remained unreturned in zone $j$ by users with profile $p$ |

$$\min_{R_{ijp}, Y_i,} \sum_{i=1}^{m}\sum_{j=1}^{n}\sum_{p=1}^{p}(x_{ijp}R_{ijp})/18 + \sum_{i=1}^{m}f_i Y_i + \sum_{j=1}^{m}\sum_{p=1}^{p}mb_{jp} \qquad (1)$$

**Subject to,**

$$w_{ijp} \leq Y_i A_{ijp}, \qquad \forall\, i \in I, \forall\, j \in J, \forall\, p \in P \qquad (2)$$

$$x_{ijp} \leq k_i w_{ijp}, \qquad \forall\, i \in I, \forall\, j \in J, \forall\, p \in P \qquad (3)$$

$$\sum_{j=1}^{n}\sum_{p=1}^{p} x_{ijp} \leq k_i Y_i, \qquad \forall\, i \in I \qquad (4)$$

$$\sum_{i=1}^{m} x_{ijp} + b_{jp} \geq \theta Q_{jp} \qquad \forall\, j \in J, \forall\, p \in P \qquad (5)$$

$$R_{ijp} \geq w_{ijp}(c_{ij} + a_p), \qquad \forall\, i \in I, \forall\, j \in J, \forall\, p \in P \qquad (6)$$

$$Y_i \in \{0,1\}, \qquad \forall\, i \in I \qquad (7)$$

$$w_{ijp} \in \{0,1\}, \qquad \forall\, i \in I, \forall\, j \in J, \forall\, p \in P \qquad (8)$$

$$x_{ijp}, R_{ijp}, b_{jp} \geq 0, \qquad \forall\, i \in I, \forall\, j \in J, \forall\, p \in P \qquad (9)$$

The above formulation presented in equations (1) – (9) resulted in a mixed-integer non-linear mathematical model. The objective function presented in equation (1) aims to minimize the total cost that consists of three cost components: (i) fixed cost for opening and operating an opioid disposal kiosk, (ii) the cost of providing incentive in accordance to the number of unused opioids per prescription that are returned, (iii) the penalty cost in the event of failing to ensure the return of unused opioids that were set as a targeted return at the beginning of the campaign. Constraint family (2) requires that users can return unused opioids in an opened disposal kiosk that is placed within the maximum distance that users will be willing to travel at given minimum reservation incentive level. Constraints (3) ensure that users from multiple zones can return unused opioids to a single disposal kiosk as long as the total number of returned opioids does not exceed the capacity of the focal disposal kiosk. Constraint sets (4) make sure that the quantity of opioids that are returned to a disposal kiosk will not exceed the capacity of that kiosk. Constraints (5) ensure sufficient return of the unused opioids commensurate with the quantity of opioids that are set as target return. Finally, constraints (6) ensure that users in a certain zone will only return opioids if the offered incentive is higher than the sum of the users' minimum reservation incentive level and travelling cost to a disposal kiosk. Finally, the binary and non-negativity nature of the decisions are imposed with the constraints (7), (8), and (9).

**4.1 Decomposition of the MINLP**

As noted, the resulting formulation given by equations (1) – (9) is a Mixed-Integer Non-Linear Programming (MINLP) model, which we implemented in AMPL and attempted to solve to optimality using KNITRO 12.4.0 solver. The computing machine has the following specifications: Dual Intel Xeon Processor (12 Core, 2.3 GHz Turbo) with 64 GB of RAM and a 64-bit operating system. We considered a smaller test case instance including 50 potential locations for placing opioid disposal kiosks and 12 zones from where the opioid users will potentially return their unused medications. However, despite running for roughly 55 hours, the model was unable to converge and provide an optimal solution, which shows computational intractability of this MINLP model for the above-mentioned smaller problem instance.

To efficiently solve this problem, we consider a reformulation of the MINLP using decomposition technique to solve the original problem in two different stages: first, and second stage. The characteristics of the decision problem are suitable for this: decisions for opioid disposal kiosk locations and incentive amounts are made in the first stage, whereas assignments of zones to kiosk locations, returned and unreturned quantities are decided in the second stage. Decisions made in the second stage understandably rely on actions taken in the first stage. This reformulation results in a model structure that is suitable to solve efficiently with Benders Decomposition (BD) approach that takes the advantage of dual decomposition structure for splitting the original MINLP in master problem (i.e., the first stage problem) and a set of sub-problems (i.e., the second stage problems). Unlike metaheuristic algorithms, Benders Decomposition (BD) approach helps us to obtain the guaranteed optimality [46]. It also allows us to use state-of-the-art solvers (CPLEX) to solve the master problem and subproblems very efficiently. The resulting decomposition scheme has the following three properties: (i) the master problem preserves the NP hardness (i.e., difficult to solve) pertinent to the Mixed-Integer Linear Programming (MILP) structure of the problem, (ii) the non-linearity feature of the problem is absorbed by the subproblems, and (iii) solving the sub-problems to the optimality ensures the zero-duality gap. We described the reformulated problem alongside the BD algorithms in the subsequent subsections.

### 4.2 Derivation of the sub-problem (2nd stage)

We formulated the sub-problem considering an initial feasible solution (i.e., the location of the disposal kiosks ($\bar{Y}_i$) and incentive to be paid for the return of unused opioids ($\bar{R}_{ijp}$) that were temporarily fixed).

$$\min_{x_{ijp},\ b_{jp}} \sum_{i=1}^{m}\sum_{j=1}^{n}\sum_{p=1}^{p}(x_{ijp}\bar{R}_{ijp})/18 + \sum_{j=1}^{m}\sum_{p=1}^{p} Mb_{jp} \tag{10}$$

**Subject to,**

$$w_{ijp} \leq \bar{Y}_i A_{ijp}, \qquad \forall\ i \in I, \forall\ j \in J, \forall\ p \in P \tag{11}$$

$$x_{ijp} \leq k_i w_{ijp}, \qquad \forall\ i \in I, \forall\ j \in J, \forall\ p \in P \tag{12}$$

$$\sum_{j=1}^{n}\sum_{p=1}^{p} x_{ijp} \leq k_i \bar{Y}_i, \qquad \forall\ i \in I \tag{13}$$

$$\sum_{i=1}^{m} x_{ijp} + b_{jp} \geq \theta Q_{jp} \qquad \forall\ j \in J, \forall\ p \in P \tag{14}$$

$$\bar{R}_{ijp} \geq w_{ijp}(c_{ij} + a_p), \qquad \forall\ i \in I, \forall\ j \in J, \forall\ p \in P \tag{15}$$

$$w_{ijp} \in \{0,1\}, \quad \forall\ i \in I, \forall\ j \in J, \forall\ p \in P \tag{16}$$

$$x_{ijp}, b_{jp} \geq 0, \quad \forall\ i \in I, \forall\ j \in J, \forall\ p \in P \tag{17}$$

The objective function (10) minimizes the total cost of providing an incentive for the return of unused opioids and also the penalty cost resulting from the quantity of unreturned opioids. The constraint family (11), (12), and (13) put the same restriction as (2), (3), and (4) did in the original MINLP model to restrict that opioid user with a certain profile in a certain zone will return unused opioids to an opened disposal kiosk (i.e., the solution of the master problem) that has sufficient capacity. Similar to (5), constraints (14) ensure sufficient return of the unused opioids. Finally, similar to (6), constraint family (15) consider the incentive determined by the master problem and make sure that the users will only return their unused opioids if that incentive is higher than the sum of the travelling cost and minimum reservation incentive level of the user. The constraints (16) and (17) restricts the binary nature and non-negativity of the decisions, respectively.

### 4.3 Derivation of the master problem (1$^{st}$ stage)

The master problem only considered the decision variables pertinent to opening the disposal kiosk at a particular location ($Y_i$) and amount of incentive ($R_{ijp}$) to be paid in return of the unused opioids. Constraints family that included only these two decision variables were included in the master problem. In addition, we introduced a new variable $u$ to iteratively build and add an optimality cut to the master problem according to BD algorithm.

Here, the objective function (18) minimizes the total cost associated with the opening of kiosk locations as well as the cost of incentive allocation and penalty of failure to return, which is represented by the artificial variable $u$. The constraints (19) update $R_{ijp}$ ensuring the same restriction as imposed by constraints (6), where $\bar{w}_{ijp}$ is the parameterized solutions coming from the 2$^{nd}$ stage model. Constraint (20) represent the optimality cut that is iteratively generated from the solution of the 2$^{nd}$ stage model and added to the 1$^{st}$ stage. Here, $\alpha_{ijp}$, $\beta_{ijp}$, $\gamma_i$, $\delta_{jp}$, $\varphi_{ijp}$ are the dual variables corresponding to equations (11), (12), (13), (14), and (15), respectively. Finally, constraints (21) and (22) restrict the binary nature of the location and non-negativity nature of the incentive amount to be paid to the users.

$$\min_{Y_i} u + \sum_{i=1}^{m} f_i Y_i \qquad (18)$$

**Subject to,**

$$R_{ijp} \geq \bar{w}_{ijp}(c_{ij} + a_p), \qquad \forall\, i \in I, \forall\, j \in J, \forall\, p \in P \qquad (19)$$

$$u \geq \sum_{i=1}^{m}\sum_{j=1}^{n}\sum_{p=1}^{p}\alpha_{ijp}Y_{i}A_{ijp} + \sum_{i=1}^{m}\sum_{j=1}^{n}\sum_{p=1}^{p}\beta_{ijp}C_{i}\overline{w}_{ijp} + \sum_{i=1}^{m}\gamma_{i}C_{i}Y_{i} + \sum_{j=1}^{m}\sum_{p=1}^{p}\delta_{jp}\theta Q_{jp}$$

$$+ \sum_{i=1}^{m}\sum_{j=1}^{n}\sum_{p=1}^{p}\varphi_{ijp}\overline{w}_{ijp}(c_{ij}+a_{p}) \qquad \forall\, i \in I, \forall\, j \in J, \forall\, p \in P \tag{20}$$

$$Y_i \in \{0,1\}, \qquad \forall\, i \in I \tag{21}$$

$$R_{ijp} \geq 0, \qquad \forall\, i \in I, \forall\, j \in J, \forall\, p \in P \tag{22}$$

In the following subsections, we presented the steps of the benders decomposition algorithm to solve the proposed reformulated model and outlined these steps in Figure 1.

**Step 1 (initializing):** We temporarily fixed the disposal kiosks' location ($\overline{Y}_i$) and incentive to be paid to the users ($\overline{R}_{ijp}$), potentially resulting a subproblem (10) – (17) that was equivalent to the original problem except the fact that this subproblem only contained decision variables pertinent to assigning users from a certain zone to certain location, and quantity of unused opioids (both returned and unreturned).

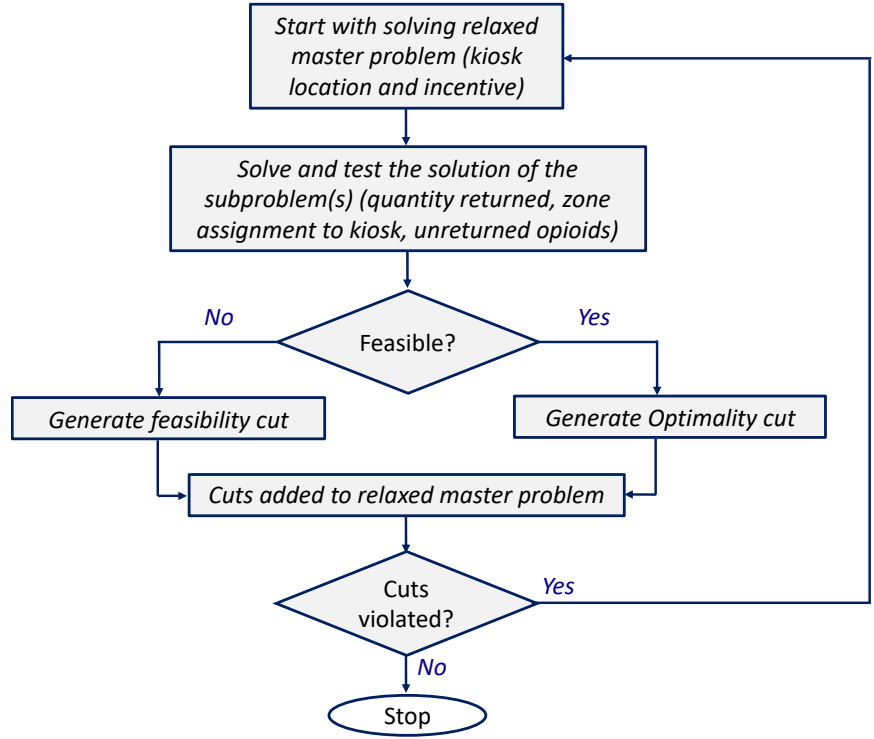

Figure 1. Outlines of Benders Decomposition algorithm

- The master problem (18) – (22) is the relaxed version of the original MINLP comprising only those variables we previously fixed (binary integer decisions for locating disposal kiosks, and the continuous variable for deciding the amount of incentive to be paid to the opioid users) and constraint family that involved those decisions.
- A set of feasibility and optimality cuts (20) are built and added to this master problem in an iterative fashion to guarantee the feasibility and optimality for the original problem.

**Step 2:** We solved the sub-problem (10) – (17) for the decision variables related to assigning users from a certain zone to a certain location, and quantity of unused opioids (both returned and unreturned). Solving the sub-problem also produce the optimal dual variables ($\alpha_{ijp}$, $\beta_{ijp}$, $\gamma_i$, $\delta_{jp}$, $\varphi_{ijp}$).

**Step 3:** We built and added a new optimality cut (20) to the master problem depending on the optimal solutions (i.e., primal and dual variables) obtained after solving the subproblem.

**Step 4:** We solved the master problem (18) – (22) with the new cuts and obtained a set of decisions for locating disposal kiosks and incentive to be paid to the opioid users.
- Case (i): These new set of decisions (locating disposal kiosks and incentive to be paid) are equal to the decisions obtained in the previous iteration. This guaranteed the feasibility and optimality conditions for the original problem, and we terminated the algorithm.
- Case (ii): If these decisions are not equal, then we went to Step 2 (i.e., solved the sub-problem once again keeping the new set of master problem's decisions fixed).

This decomposition scheme ensured convergence of the master problem, which in turn guaranteed zero-duality gap for the original problem. Specifically, the master problem solution is considered as lower bound of the original problem, which together with the global optimal conditions of the sub-problem in each iteration ensured the zero-duality gap.

## 5. Numerical experiment and discussion

In this section, we present a test case for our proposed MINLP model (1) – (9) and show the results after solving this test case problem using BD algorithm. We implemented the reformulated model (10) –

(22) in AMPL and solved it to the optimality using CPLEX 20.1.0.0 solver in the same computer that we used to solve the proposed MINLP model (1) – (9). We constructed this test case problem for the Massachusetts Middlesex county which experienced the highest number of opioid-related overdose deaths in the state over the period of last ten years (i.e., from 2010 to 2020). This county has a total of 138 existing pharmacy (CVS and Walgreens) locations that we used as potential sites for opening opioid-disposal kiosks. The annual fixed cost for opening a disposal kiosk with a capacity of 30,000 opioids (i.e., number of pills) is assumed to be $2000 [40, 47]. The opioid users in this county are segregated in twelve major cities that were considered as different zones. We first computed the travelling distance from the centroids of each zone to each pharmacy location. A value of 50 cents per mile of travelling distance were assumed as cost of transportation to estimate the travelling cost from a certain zone to a certain location [48]. A portion of this travelling cost is presented in Table 2. We considered a penalty cost of $12 for each prescription for which the unused opioids were not returned.

| Table 2. Cost of travelling from a certain location to certain zones | | | | | | | | |
|---|---|---|---|---|---|---|---|---|
| **Pharmacy locations** | **Cambridge** | **Everett** | **Framingham** | **Lowell** | **Malden** | **Marlboro** | **Medford** | **Melrose** |
| 344 Great Road Acton, MA 01720 | 16.4 | 25.2 | 21.75 | 13.15 | 23.9 | 19.25 | 15.95 | 21.7 |
| 400 Massachusetts Ave, Acton, MA 01720 | 16.3 | 25.1 | 12.55 | 15.9 | 23.8 | 18.1 | 15.85 | 21.6 |
| 23-25 Massachusetts Ave, Alwife Plz Arlington, MA 02474 | 2.8 | 4.35 | 18.85 | 20 | 4.65 | 26.4 | 2.4 | 6.65 |
| 833 Massachusetts Avenue, Arlington, MA 02476 | 4.45 | 5.45 | 18.1 | 19.2 | 4.35 | 25.65 | 2.55 | 6.2 |
| 47 Pond St, Kings Crossing Plaza Ashland, MA 01721 | 16.9 | 22.55 | 1.4 | 29.2 | 24 | 9.7 | 23.5 | 27.3 |
| 414 Union St. Ashland, MA 01721 | 20.95 | 26.55 | 2.8 | 32.05 | 28 | 7.15 | 27.55 | 31.3 |
| 199 Great Road Bedford, MA 01730 | 12.1 | 15.3 | 17.15 | 10.75 | 14 | 24.7 | 12.55 | 11.8 |
| 264 Trapelo Rd. Belmont, MA 02478 | 3.35 | 7.75 | 13.3 | 19.35 | 6.7 | 20.85 | 4.9 | 8.55 |
| 60 Leonard Street Belmont, MA 02478 | 3.45 | 6.4 | 17.45 | 18.55 | 5.3 | 24.95 | 3.5 | 7.15 |
| 700 Boston Rd, Rt 3a, Towne Plz Billerica, MA 01821 | 16.9 | 14.4 | 21.9 | 8.25 | 13.1 | 29.15 | 11.65 | 10.9 |
| 34 Cambridge Street, Space 160 Burlington, MA 01803 | 13.65 | 11.3 | 18.7 | 12.45 | 10.05 | 26.25 | 8.6 | 7.85 |
| 242 Cambridge St. Burlington, MA 01803 | 14.55 | 12 | 19.55 | 11.5 | 10.7 | 27.1 | 9.3 | 8.5 |
| 624 Massachusetts Ave. Cambridge, MA 02139 | 0.9 | 3.4 | 15.2 | 23 | 4.4 | 22.7 | 4 | 7.35 |
| 215 Alewife Brook Parkway Cambridge, MA 02138 | 1.95 | 5.25 | 19.2 | 20.3 | 5.55 | 26.75 | 3.3 | 7.6 |
| 225 Cambridge St Cambridge, MA 02141 | 1.25 | 2.85 | 17.3 | 24.4 | 4 | 24.85 | 4.65 | 8 |
| 100 Cambridgeside Place, Suite E122 Cambridge, MA 02141 | 1.6 | 2.8 | 17 | 24.35 | 3.8 | 24.55 | 4.6 | 7.95 |
| 6 Jfk Street Cambridge, MA 02138 | 0.65 | 3.8 | 15.3 | 21.65 | 4.75 | 22.85 | 3.15 | 11.35 |
| 36 White Street Cambridge, MA 02140 | 1.25 | 4 | 16.2 | 20.8 | 4.65 | 23.75 | 2.3 | 7.15 |
| 16 Boston Road Chelmsford, MA 01824 | 21.45 | 22.5 | 26.5 | 3.85 | 21.25 | 22.05 | 19.8 | 19.05 |
| 199 Sudbury Road Concord, MA 01742 | 12.6 | 21.4 | 18 | 13.05 | 20.15 | 21.9 | 12.15 | 17.95 |

For each zone, we computed the number of opioid prescriptions in accordance with the twelve profiles of opioid users. We used MA APCD to extract the claims related to opioid prescriptions for the year 2014. In order to estimate the number of unused opioids per prescription, we assumed that on an average 60% of prescribed opioids per prescription remained unused—an anecdotal proportion that was reported in the previous study [23]. We present the quantity of unused opioids that we estimated for different zones in accordance to different user profiles in Table 3.

| Table 3. Unused opioids potentially available for return in different zones | | | | | | | | | | | | |
|---|---|---|---|---|---|---|---|---|---|---|---|---|
| Opioid users' zones | Profile 1 | Profile 2 | Profile 3 | Profile 4 | Profile 5 | Profile 6 | Profile 7 | Profile 8 | Profile 9 | Profile 10 | Profile 11 | Profile 12 |
| Cambridge | 9198 | 82560 | 21606 | 193866 | 11328 | 101100 | 5556 | 51876 | 16644 | 151164 | 6828 | 59616 |
| Everett | 4950 | 46164 | 17892 | 161598 | 5832 | 53160 | 3360 | 29208 | 15228 | 141810 | 3306 | 30528 |
| Framingham | 6768 | 63972 | 23214 | 214542 | 12336 | 111456 | 4926 | 43752 | 18906 | 171444 | 6528 | 59712 |
| Lowell | 19188 | 173082 | 53700 | 489672 | 21582 | 190050 | 11028 | 99384 | 43200 | 391104 | 11412 | 101784 |
| Malden | 6918 | 62784 | 23484 | 210750 | 8892 | 80352 | 5058 | 44178 | 19518 | 176898 | 5286 | 47448 |
| Marlboro | 5256 | 46788 | 17010 | 155244 | 8436 | 72984 | 3876 | 36300 | 13596 | 116226 | 4086 | 37824 |
| Medford | 8556 | 77454 | 21858 | 200844 | 11748 | 106476 | 6066 | 50664 | 21606 | 195480 | 6060 | 57564 |
| Melrose | 2784 | 24024 | 9480 | 83706 | 5556 | 49026 | 1896 | 17670 | 8028 | 72864 | 3012 | 26358 |
| Newton | 4044 | 40230 | 18624 | 162810 | 10836 | 96714 | 3576 | 31902 | 14676 | 128268 | 6474 | 59286 |
| Somerville | 8562 | 78930 | 21432 | 188130 | 8826 | 81966 | 6348 | 56946 | 17688 | 160542 | 5040 | 46164 |
| Waltham | 6846 | 62520 | 21702 | 191532 | 9522 | 86586 | 5472 | 47262 | 16374 | 145836 | 5766 | 50862 |
| Woburn | 5712 | 50988 | 17178 | 149736 | 7476 | 66924 | 4020 | 37824 | 14484 | 126000 | 5166 | 45450 |

As noted, we particularly included user profiles to take into account the variation in users' willingness to return unused opioids. These profiles were partly informed by our previous study where we reported that patients' age and gender were important predictor of opioid use disorder [45]. We also observed that individuals who are male and have age in between 18 to 35 were more likely to develop opioid use disorder (OUD) and related outcomes. OUD is caused by the abuse and misuse of opioids, which often are associated with the diversion of prescription opioids. Thus, we assumed that individuals who may potentially be at an elevated risk of developing OUD will also be less likely to return their unused opioids. We classified age into three categories: 18 to 35, 36 to 45, and above 45 to capture the variation in users' motivation for returning unused opioids in different age groups. Previous study also indicated that individual who were prescribed immediate-release opioid were twice as unlikely to return their unused opioids compared to those who were prescribed extended-release opioids [1]. Therefore, we considered individuals' gender, age, and type of prescribed opioids to construct their profiles. The combinations of these factors resulted in twelve profiles of opioid users. For each profile, the estimated three levels of minimum reservation incentive (i.e., low, medium, and high incentive) are presented in Table 4. We considered a maximum distance of 4, 8, and 20 miles that an opioid user will be willing to travel for a low, medium, and high reservation incentive, respectively. Finally, we solved the test case problem for three different values (i.e., 50%, 80%, and 100%) of $\theta$ (i.e., the percentage of potentially available unused opioids that the third-party company intended to return through the campaign).

| Table 4. Profile specific minimum reservation incentive levels ||||
|---|---|---|---|
| **Profiles of opioid users** | **Low minimum Reservation Incentive** | **Medium minimum Reservation Incentive** | **High minimum Reservation Incentive** |
| Profile 1 | 10.5 | 12.5 | 14.5 |
| Profile 2 | 13.5 | 15.5 | 17.5 |
| Profile 3 | 7.5 | 9.5 | 11.5 |
| Profile 4 | 9 | 11 | 13 |
| Profile 5 | 4.5 | 6.5 | 8.5 |
| Profile 6 | 6.5 | 8.5 | 10.5 |
| Profile 7 | 12.5 | 14.5 | 16.5 |
| Profile 8 | 15 | 17 | 19 |
| Profile 9 | 9.5 | 11.5 | 13.5 |
| Profile 10 | 11 | 13 | 15 |
| Profile 11 | 6.5 | 8.5 | 10.5 |
| Profile 12 | 8 | 10 | 12 |

For each value of $\theta$, we present the annual cost for running the campaign in Figure 2 (a). A breakdown of this total cost into the fixed cost for opening optimal number of disposal kiosks, cost of incentive, and penalty cost associated with the unreturned opioids is presented in Figure 2 (b), 2 (c), and 2 (d), respectively. We observed that the fixed kiosk-opening cost and incentive cost show an increasing trend with the increase in minimum reservation incentive level for each value of $\theta$. However, we observed that the total cost for running the campaign and penalty cost attributable to the quantity of unreturned opioids decreased with an increase in the minimum reservation incentive level. These findings indicate that providing a higher incentive is cost effective to run the overall campaign. Further investigation explains that with a higher incentive level, since the opioid users are willing to travel further, the model opened disposal kiosks at some other locations, potentially increasing the overall capacity of the disposal kiosks. Although it increases the total amount of incentive that are paid to the users due to the increased number opioids that can be returned at certain value of $\theta$, it decreases the penalty cost associated with the unreturned opioids, which contribute to the reduction in the total cost.

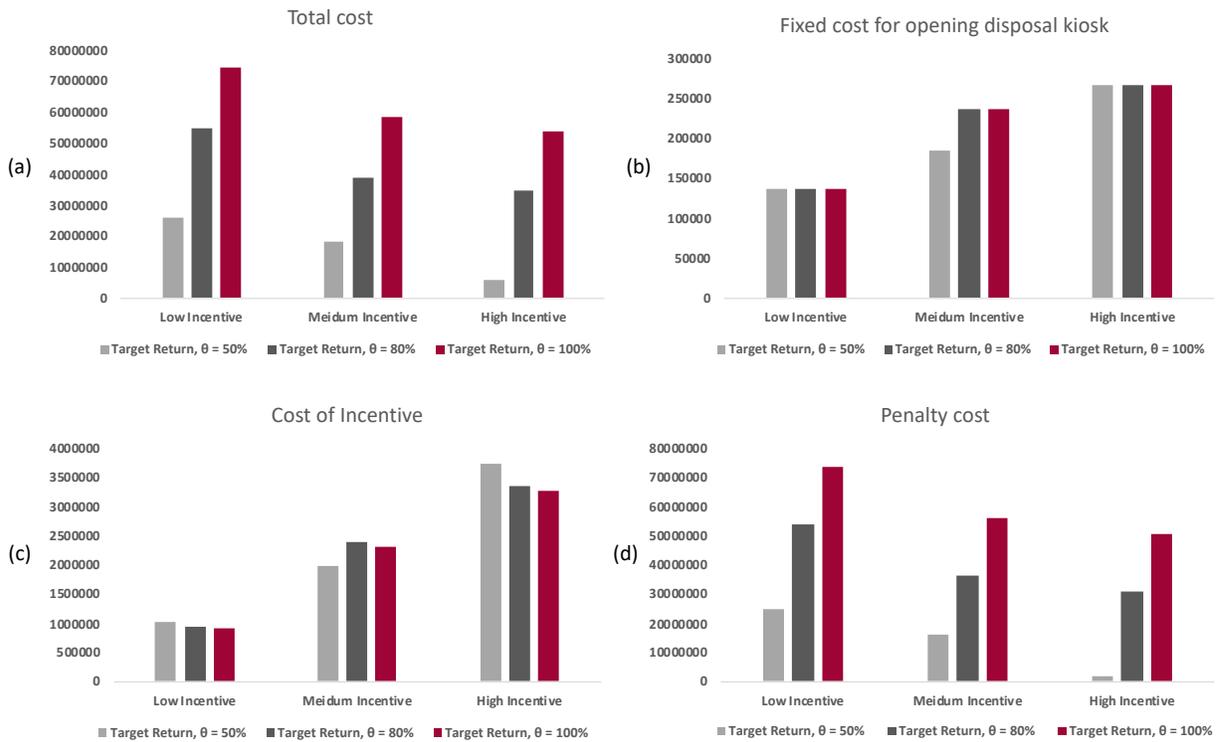

Figure 2. Breakdown of different cost components to run the campaign. (a) shows the total cost (upper left), (b) shows the fixed cost for opening a disposal kiosk (upper right), (c) shows cost of incentive (lower left), and (d) shows the penalty cost (lower right)

For a target return quantity ($\theta$) of 50%, we also graphically presented the optimal locations of opioid disposal kiosks where users with different profiles returned their unused opioids. For better representation, we only displayed the disposal kiosks' locations that were selected for two zones: Cambridge and Somerville. The kiosks assignment for users from Cambridge and Somerville are presented using red and black lines, respectively. We showed the optimal kiosks' locations in these two cities when users were offered with the low and high reservation incentives in Figure 3 and Figure 4, respectively. The idea is to demonstrate the influence of providing higher incentives on users' willingness to travel further to return the opioids. We refer here to the assumption that users were willing to travel less distance at the low incentive level compared to what they were willing to travel at the high incentive level. In compliance to this assumption, we observed that a fewer number of kiosks were placed in each of these two cities when offered with the low incentives (Figure 3). Since a fewer number of disposal kiosks are placed, a large number of unused opioids remained unreturned due to the less kiosks' capacity, leading to a high penalty cost. However, when offered with the high incentives, we observed that a few more disposal kiosks are placed at a far distance compared to the locations selected at the low incentives (Figure 4). This increased the overall capacity of the disposal kiosks and reduced the quantity of unused opioids that remained unreturned at the low incentive level, which in turn resulted in the lower penalty cost. We also observed that regardless of the offered incentive levels, users with different profiles from two neighborhood cities (i.e., Cambridge and Somerville) shared a single disposal kiosk for returning their unused opioids.

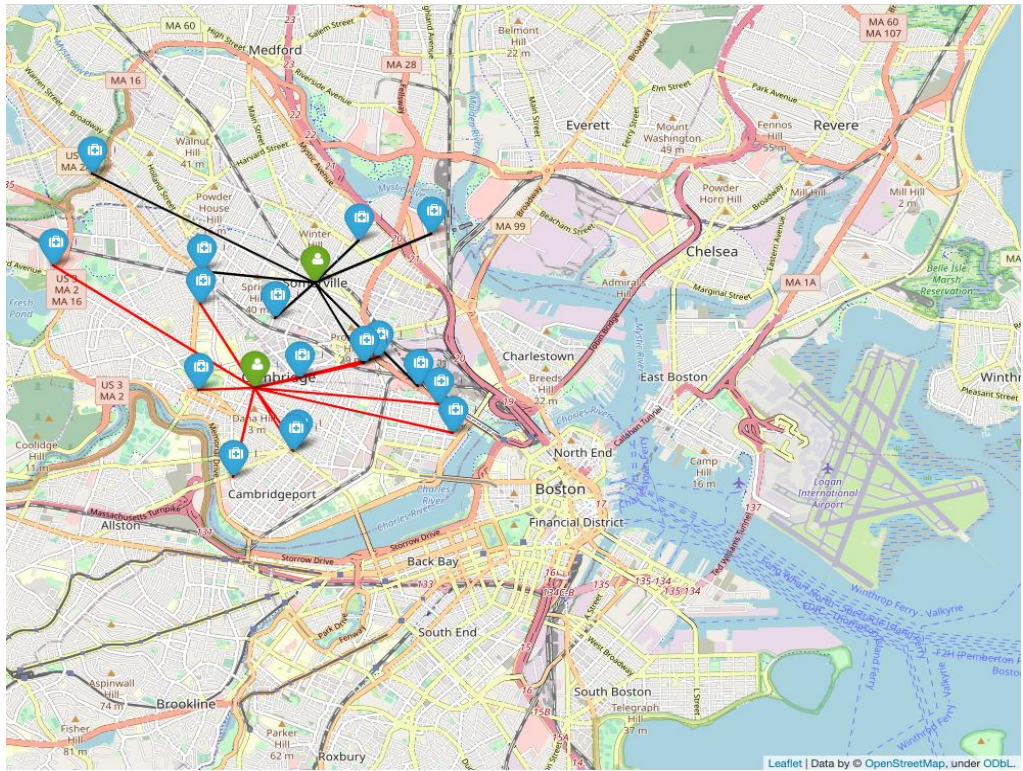

Figure 3. Optimal kiosks' location at the low reservation

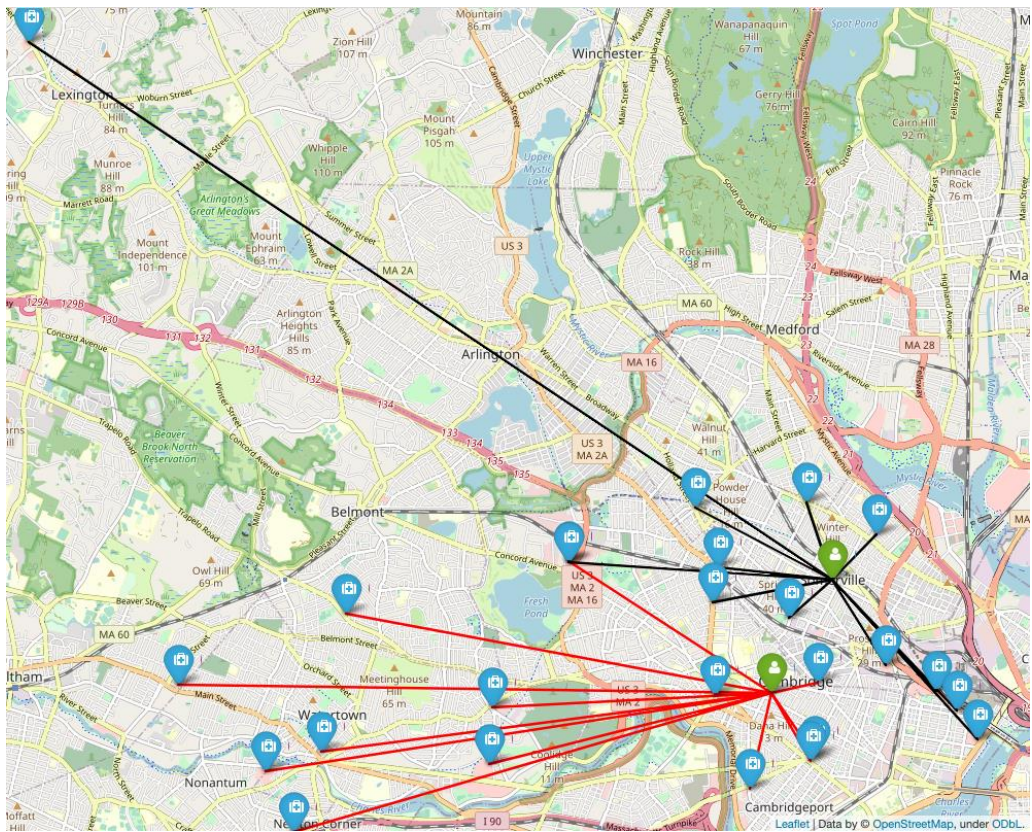

Figure 4. Optimal kiosks' location at the high reservation incentive

## 6. Summary, limitations, and outlook

In this study, we aimed to address the problem pertaining to the diversion of prescription opioids to secondary users and the black market. Our proposed framework was to locate opioid disposal kiosks at sites accessible to prescription opioid users and to incentivize such users to return their unused prescription opioids. To that end, we formulated a mathematical model that can help decide the locations of disposal kiosks and the amount of incentive to be paid to the prescription opioid users for optimizing the return and secure disposal of prescription opioids. Given the fact that the proposed MINLP model was computationally challenging to solve, we reformulated the problem using BD technique. We solved this reformulated model by applying it to the Massachusetts Middlesex county using data from the MP APCD. This numerical experiment showed the model's ability in providing insights that can potentially help state public health officials design various trade-off strategies for implementing a feasible and sustainable campaign for returning unused opioids. Such a program, if implemented successfully, can potentially reduce the diversion of prescription opioids to the drug traffickers and prevent the non-medical use of legitimate opioids by secondary users.

Despite this contribution, we acknowledge several limitations of this data-driven optimization framework, which also provides directions for extending this research. First, although some parameters (e.g., patient profiles and profile-specific parameters, and estimation of potentially unused opioids) could be transferable to non-US healthcare settings, the direct implementation of our proposed optimization framework in non-US healthcare settings may lack generalizability. We designed our problem based on the pharmaceutical stewardship laws passed in several US states that require pharmaceutical manufacturers to facilitate and fund drug stewardship programs for collecting prescription and non-prescription drugs [40, 41]. The penalty cost (due to the stewardship organization's failure for ensuring the return of unused opioids) that we considered in our model is also based on that drug stewardship law [40, 41]. The possibility of placing kiosks in a retail pharmacy (such as CVS and Walgreens) is also linked to the US context [41]. Therefore, to adapt our model to other healthcare systems, it would be necessary to tailor it to the specific roles played by key stakeholders (i.e., government, state healthcare agency, and pharmaceutical company, and retail pharmacy) in that healthcare system. While our goal in this study was to present a conceptual framework that can be used to reduce the diversion of prescription opioids, the implementation of this framework/model to other healthcare systems (other than US) will require further refinements of several parameters used in our proposed model.

Second, we lacked data pertaining to setting the minimum reservation incentive and the maximum distance users will be willing to travel to return their unused drugs. To set a reasonable estimate of these two parameters, we largely followed the findings reported in a previous study [23, 24] and assumed the values by adjusting the reported incentive amount based on current economy. An alternative approach would be fitting a parametric distribution based on some assumptions to obtain an approximation of this minimum reservation incentive. However, the fact that users' motivation to return their unused drugs is influenced by different socio-behavioral characteristics, assuming parametric distributions in this regard may not be realistic. Because these two parameters are crucial to analyze the distance and incentive tradeoffs, which in turn influence the overall cost estimation for administering the campaign, further research can be conducted to collect actual data regarding the two parameters through a pilot study. Future research can also include different sources of uncertainty in parameter estimation. Such uncertain parameters can be the number of unused opioids that may potentially be available in a certain zone, the minimum reservation incentive level, and the maximum distance that an opioid user with a certain profile will be willing to travel to return their unused opioids. Multi-stage stochastic programming methods can be adopted to model such uncertainty and formulate a more robust and pragmatic decision-making framework.

To summarize, we expect that our proposed optimization framework and study findings will help strengthen community resilience against the opioid addiction epidemic in conjunction with other efforts, contributing to the health and well-being of individuals, communities and societies. By leveraging the cross-scaler engagement with local and national stakeholders and public health policy makers, our study can inform implementable policies to reduce the diversion of unused prescription opioids. With this unique decision-making framework, we hope to expand the existing knowledge base for addressing complex dynamics of the opioid addiction epidemic.

36. Shaw, P. and S. Maynard, *The potential of financial incentives to enhance householders' kerbside recycling behaviour.* Waste Management, 2008. **28**(10): p. 1732-1741.
37. Fabbricino, M., *An integrated programme for municipal solid waste management.* Waste management & research, 2001. **19**(5): p. 368-379.
38. Weraikat, D., M.K. Zanjani, and N. Lehoux, *Two-echelon pharmaceutical reverse supply chain coordination with customers incentives.* International Journal of Production Economics, 2016. **176**: p. 41-52.
39. Taleizadeh, A.A., E. Haji-Sami, and M. Noori-daryan, *A robust optimization model for coordinating pharmaceutical reverse supply chains under return strategies.* Annals of Operations Research, 2019: p. 1-22.
40. Cuomo, A.M. and B. Seggos, *Statewide Pharmaceutical Stewardship Program.*
41. Massachusetts Department of Public Health (MDPH). *DRUG STEWARDSHIP PROGRAM*. 2021 [cited 2021 June 03]; Available from: https://www.mass.gov/doc/inmar-rx-solutions-proposed-drug-stewardship-plan/download.
42. The Commonwealth of Massachusetts. *DRUG STEWARDSHIP PROGRAM*. 2016; [cited 2022 September 21] Available from: https://www.adea.org/uploadedFiles/ADEA/Blogs/ADEA_State_Update/MA_H4056.pdf.
43. de Figueiredo, J.N. and S.F. Mayerle, *Designing minimum-cost recycling collection networks with required throughput.* Transportation Research Part E: Logistics and Transportation Review, 2008. **44**(5): p. 731-752.
44. Dufour, R., et al., *Understanding predictors of opioid abuse: predictive model development and validation.* Am J Pharm Benefits, 2014. **6**(5): p. 208-216.
45. Hasan, M.M., et al., *A machine learning framework to predict the risk of opioid use disorder.* Machine Learning with Applications, 2021. **6**: p. 100144.
46. Blum, C. and A. Roli, *Metaheuristics in combinatorial optimization: Overview and conceptual comparison.* ACM computing surveys (CSUR), 2003. **35**(3): p. 268-308.
47. SHARPS. Compliance Inc. *MedSafe*. 2016 [cited 2021 June 5]; Available from: https://ncweb.pire.org/documents/Temporary%20ECHO%20sites/MedSafe%20Gov_StateLocal%20Pricing.pdf.
48. Internal Revenue Service (IRS). *IRS issues standard mileage rates for 2021*. 2021 [cited 2021 May 12]; Available from: https://www.irs.gov/newsroom/irs-issues-standard-mileage-rates-for-2021.
28